\documentclass{article}
\usepackage{amssymb}
\usepackage{amsmath}
\usepackage[pdftex]{color}
\newtheorem{thm}{Theorem}[section]
\newtheorem{lem}[thm]{Lemma}
\newtheorem{prop}[thm]{Proposition}
\newtheorem{cor}[thm]{Corollary}
\newcommand{\C}{\mathbb C}
\newcommand{\Z}{\mathbb Z}
\newcommand{\N}{\mathbb N}
\newcommand{\hm}{{\rm Hom}}
\newcommand{\gr}{{\rm gr}}
\newcommand{\der}{{\rm Der}}
\newcommand{\ap}{{\rm Ap}}
\newcommand{\val}{{\rm val}}
\newcommand{\R}{\mathbb C[S]}
\newcommand{\A}{\mathbb C[T]}
\title{On differential operators of numerical semigroup rings}
\author{Valentina Barucci, Ralf Fr\"oberg}
\begin{document}
\maketitle
\begin{abstract}
If $S=\langle d_1,\ldots,d_\nu\rangle$ is a numerical semigroup, we call
the ring $\C[S]=\C[t^{d_1},\ldots,t^{d_\nu}]$ the semigroup ring of $S$.
We study the ring of differential operators on $\C[S]$, and its
associated graded in the filtration induced by the order of the differential
operators. We find that these are easy to describe in case $S$ is a so called
Arf semigroup. If $I$ is an ideal in $\C[S]$ that is generated by monomials,
we also give some results on $\der(I,I)$ (the set of derivations which map $I$
into $I$).
\end{abstract}

 \section{Introduction}
It is well known that the ring $D(\R)$ of differential operators of a semigroup ring $ \R=\C[t^{d_1},\ldots,t^{d_\nu}]$, where  $\C$ is the complex field,  is a subring of $D( \C[t,t^{-1}])$, the ring of   differential operators of the Laurent polynomials $ \C[t,t^{-1}] $. $D( \C[t,t^{-1}])$ is a non commutative $\C[t,t^{-1}]$-algebra generated by    $\partial$, the usual derivation $d/dt$. It is also known that $D(\R)$ is a $\C$-algebra finitely generated and a set of generators was found  independentely   in \cite{EE} and \cite{AE}. 
 The ring  $D(\R)$ inherits a grading from
the ring of differential operators of $ \C[t,t^{-1}]$, where $\deg(t^s)=s$ and $\deg(\partial)=-1$.
Its associated graded is a commutative
Noetherian   subring of the ring of polynomials in two indeterminates $\C[t,y]$ and it is a semigroup ring $\C[\Sigma]$, where $\Sigma \subseteq \N^2$ is a semigroup, with $|\N^2 \setminus \Sigma|$ finite. It follows that $D(\R)$  is right and left Noetherian.  
 
 In this paper we study that commutative associated ring $\C[\Sigma]$, in terms of the starting numerical semigroup $S$. Many properties of the semigroup $\Sigma$, including the minimal set of generators, can be predicted looking at $S$. If $S$ is of maximal embedding dimension, then $\Sigma$ behaves well with respect to the blowup of the maximal ideal. If,  moreover, $S$ is an Arf semigroup, then  we show how $\Sigma$, and the ring $\C[\Sigma]$ as well, is completely determined by $S$. In Section 4 we characterize the irreducible ideals of $\Sigma$, i.e. the irreducible monomial ideals of $\C[\Sigma]$ and determine the number of components for a principal monomial ideal as irredundant intersection of irreducible ideals. Finally in Section 5 we observe that, if $I$ is a monomial ideal of $\R$, then $\der(I,I)$, the $\R$-module of derivations which map $I$ into $I$ is isomorphic to the overring $I:I$ of $\R$. Thus we study the overrings of this form in relation with the monomial ideals $I$ which realize them.

\section{Numerical semigroups}
We fix for all the paper the following notation. $S$ is a  {\it numerical 
semigroup}, i.e. a subsemigroup of $\mathbb N$, with zero and with finite 
complement $H(S)=\N\setminus S$ in $\mathbb N$. The numerical semigroup 
generated by 
$d_1, \dots,d_\nu \in \mathbb N$ is $S=\langle d_1, \dots,d_\nu \rangle = 
\{ \sum_{i=1}^\nu n_id_i; n_i \in \mathbb N \}$. $M=S \setminus \{0\}$ is the  
{\it maximal ideal} of $S$, $e=e(S)$ is the  {\it multiplicity} of $S$, that 
is the smallest positive integer of $S$, $g=g(S)$ is the  
{\it Frobenius number} 
of $S$, that is the greatest integer which does not belong to $S$, $n=n(S)$ is 
the number of elements of $S$ smaller of $g$. Thus we have $S=\{ 0=s_0<s_1=
e<s_2\dots<s_{n-1}<s_n=g+1,g+2, \dots\}$.  

 A {\it relative ideal} of $S$ is a nonempty subset $I$ of $\mathbb
Z$ (which is the quotient group of $S$) such that $I+S \subseteq I$
and $I+s \subseteq S$, for some $s \in S$. A relative ideal which is
contained in $S$ is an {\it integral ideal} of $S$.

 If $I$, $J$ are relative ideals of $S$, then the following is a relative 
ideal too:
 $I-J=\{ z \in \mathbb Z \ | \ z+J \subseteq I\}$.

If $I$ is a relative ideal of $S$, then $I-I$ is the biggest semigroup $T$ 
such that $I$ is an integral ideal of $T$. There is a chain of semigroups 
$S \subseteq (I-I) \subseteq (2I-2I) \subseteq \dots \subseteq \mathbb N$, 
which stabilizes on a semigroup ${\cal B} (I)=( hI-hI)$ for $h>>0$, called 
the {\it blowup} of $I$.
Setting $S_1= {\cal B}(M)$ and $S_{i+1}={\cal B}(M_i)$, where $M_i$ is the maximal 
ideal of $S_i$, the multiplicity sequence of $S$ is $(e_0,e_1,e_2,\ldots)$, 
where $e_i=e(S_i)$. 
If $z \in \mathbb Z$, set $S(z)=\{ s\in S;s\ge z\}$, which is an ideal of $S$. 
If, with the notation above, $x=s_i$, we denote for simplicity $S(s_i)$ with 
$I_i$.  

\medskip\noindent
{\bf Definitions} For each $z\in\Z$, we define the valency of $z$ with respect 
to a semigroup $S$ as val$_S(z)=|\{ s\in S; z+s\notin S\}|$.  
Let $V_i(S)=\{ a\in{\Z}; {\rm val}_S(a)\le i\}$. When there is no ambiguity about the semigroup $S$, we will write simply val$(z)$ and $V_i$ respectively.

\begin{lem}\label{base}
{\rm (a)} $S=V_0$.

{\rm (b)} $S-M=V_1$.

{\rm (c)} For each $z \in \mathbb Z$, ${\rm val}(-z)={\rm val}(z)+z$.

{\rm (d)} $I_i-I_i \subseteq V_i$.

{\rm (e)} If $a \in \N$, then ${\rm val}(a) \le n$.

 {\rm (f)} $V_i$  is a relative ideal of $S$.

\end{lem}

\noindent
 {\bf Proof.} (a) and (b) follow immediately from the definitions.  

\noindent
(c) is proved in \cite[Lemma 3.3]{EE} and   \cite{AE}.

\noindent
(d): Let 
$a \in I_i-I_i$. Then 
$a+s_j\in I_i \subseteq S$ if $j\ge i$, so at most
$\{ a+s_0,a+s_1,\ldots,a+s_{i-1}\}$ are not in $S$. 

\noindent
(e): By (d), $I_n-I_n =\N \subseteq V_n$, so  ${\rm val}(a) \le n$ for  each $a \in \N$.

\noindent
(f): Assume $a \in V_i$ 
and $s \in S$, then val$(a+s) \leq i$. In fact if $(a+s)+b \notin S$, with 
$b \in S$, also $a+(s+b) \notin S$, with $s+b \in S$. Moreover by (c) $V_i$ 
has a minumum $m \leq 0$. Thus $s+V_i \subseteq S$, for some $s \in S$. Take 
for example $s=g+1-m.$  $\Box$

\medskip

We denote by $T(S)$ the set $\{ x \in \Z ; x \notin S, x+M \subseteq M\}$, called in \cite{ROGA} the set of pseudo-Frobenius numbers. With the notation above, $T(S)=V_1 \setminus V_0$. The cardinality of $T(S)$ is the {\it type} $t$ of the semigroup $S$. It is well known that $t \le e-1$ and $t=e-1$ if and only if 
 the   numerical semigroup $S$ is of maximal embedding dimension, i.e. when the number $\nu$ of generators equals the multiplicity $e$, i.e.  when $|M\setminus 2M|=e$,
cf. e.g. \cite[Corollaries 2.23 and 3.2, 3)]{ROGA}. It is also well known that $S$ is of maximal embedding dimension if and only if  $ M-e$ is a semigroup (cf. e.g. \cite[Proposition 3.12]{ROGA}).  In this case we have $S_1=M-e$. Setting $\pm H(S)=\{\pm h; h \in \N \setminus S\}$, we have

\begin{lem}\label{MED} Let $S$ be a numerical semigroup of maximal embedding dimension. Then val$_{S_1}(z)=\text{val}_S(z)-1$, for each $z \in \pm H(S)$.

\end{lem}

\noindent
 {\bf Proof.} Let $z \in \pm H(S)$.
Denote by $\Omega_S(z)$ (respectively  $\Omega_{S_1}(z)$) the set of pairs $(s,s+z)$, with $s \in S$ and $s+z \notin S$ (resp. the set of pairs $(s_1,s_1+z)$, with $s_1 \in S_1$ and $s_1+z \notin S_1$). By definition of valency, val$_S(z)$ is the cardinality of $\Omega_S(z)$ and val$_{S_1}(z)$ is the cardinality of $\Omega_{S_1}(z)$. If $(s,s+z)\in \Omega_S(z)$, with $s \neq 0$, then $(s-e,s-e+z)\in \Omega_{S_1}(z)$ and, conversely, if $(s_1,s_1+z)\in \Omega_{S_1}(z)$,   then $(s_1+e,s_1+e+z)\in \Omega_{S}(z)$. Thus there is a 1-1 correspondence between $\Omega_{S}(z)\setminus \{(0,z)\}$ and $\Omega_{S_1}(z)$ and the conclusion follows. $\Box$

\medskip

An {\it Arf semigroup} is a numerical semigroup
$S$ such that $I_i-s_i$ is a semigroup for each $i \geq 0$. Thus, if $S$ is Arf, then $ M-e$ is a semigroup 
and $S$ is of maximal embedding dimension.
Given an Arf
semigroup $S=\{ 0=s_0<s_1<s_2\cdots\}$, the multiplicity sequence of $S$ is
$\{ s_1-s_0,s_2-s_1,s_3-s_2,\ldots\}$.  It follows 
that the multiplicity sequence $e_0,e_1,\ldots$ of an Arf semigroup is such
that for all $i$, $e_i=\sum_{h=1}^ke_{i+h}$, for some $k\ge1$. Conversely,
any sequence of natural numbers $e_0,e_1,\ldots$ such that $e_n=1$, for 
$n>>0$, and, for all $i$,
$e_i=\sum_{h=1}^ke_{i+h}$, for some $k\ge1$, is the multiplicity sequence of
an Arf semigroup $S=\{ 0,e_0,e_0+e_1,e_0+e_1+e_2,\ldots\}$, \cite{BDF}.

\begin{lem}\label{Arf} If $S$ is an Arf semigroup, then $I_i-I_i=V_i \cap 
\mathbb N$.

\end{lem}

\noindent
 {\bf Proof.} By Lemma \ref{base} (d), $I_i-I_i \subseteq V_i \cap \mathbb N$. 
Conversely let $a \in \mathbb N$ and observe that, since $S$ is Arf, if 
$a+s_i \in S$, then also  $a+s_j \in S$ for each $j \geq i$. In fact 
$a+s_i\in S$ if and only if $a+s_i\in I_i$ if and only
if $a\in I_i-s_i$. If $ j \geq i$, since $a$, $s_j-s_i  \in  I_i-s_i$, which 
is a semigroup, then $a+s_j-s_i \in I_i-s_i$, so that $a+s_j \in S$. Thus
${\rm val}(a)=i$ if and only if $a+s_j\in S$ for $j\ge i.$ $\Box$

\section{Differential operators on  numerical semigroup rings}
Let $R$ be a commutative $\C$-algebra. The ring of differential operators
$D(R)$ of $R$ is inductively defined in the following way. Setting
$$D^0(R)=\{ \Theta_a;a\in R\}$$
$$D^i(R)=\{\Theta\in\hm_k(R,R);[\Theta,D^0(R)]\subseteq D^{i-1}(R)\}$$

\noindent where $\Theta_a\colon R\rightarrow R$ is the multiplication map $r\mapsto ar$,
and $[\Theta,\Phi]=\Theta\Phi-\Phi\Theta$ is the commutator, the ring of differential operators is 
$$D(R)=\cup_{i\ge0}D^i(R)$$
This is a filtered ring, $D^i(R)D^j(R)\subseteq D^{i+j}(R)$ for all $i$ and 
$j$, and $D^i(R)\subseteq D^{i+1}(R)$, and its associated graded is
$\gr(D(R))=\oplus_{i\ge0}D^i(R)/D^{i-1}(R)$, where $D^{-1}(R)=0$. 

\smallskip
The module of derivations $\der(R)$ is $\{\Theta\in\hm_k(R,R);\Theta(ab)=
\Theta(a)b+a\Theta(b),\ a,b\in R\}$\ and it is well known that $\der(R)$ is $\{\Theta\in D^1(R); \Theta(1)=0\}$.

\bigskip
The rings of differential operators of semigroup rings have been studied by
e.g. Perkins (\cite{PE}),   Eriksen (\cite{EE}) and Eriksson (\cite{AE}). It is shown that the ring of 
differential operators of a semigroup ring $\C[S]$ is a subring of $D(\C[t,t^{-1}])=
\C[t,t^{-1}]\langle\partial\rangle$. The ring of 
differential operators of a semigroup ring inherits a grading from
$D(\C[t,t^{-1}])$, where $\deg(t^s)=s$ and $\deg(\partial)=-1$.
Its associated graded is a commutative
Noetherian   subring of the ring of polynomials $\C[t,y]$ and it is a semigroup ring $\C[\Sigma]$, where $\Sigma \subseteq \N^2$ is a semigroup, with $|\N^2 \setminus \Sigma|$ finite.  More
precisely, it was independentely proved in \cite{EE} and \cite{AE}  that 

\begin{thm}\label{E}{\rm \cite{EE},\cite{AE} }If $S=\langle d_1,\ldots,d_{\nu}\rangle$ 
is a numerical semigroup 
and $\C[S]$ its semigroup ring, then $\gr(D(\C[S]))$ 
is a $\C$-subalgebra of $\C[t,y]$ minimally generated by the monomials 
$$\{ t^{d_1},\ldots,t^{d_{\nu}},y^{d_1},\ldots,y^{d_{\nu}}\}\cup\{ ty\}\cup
\{ t^{{\rm val}(-h)}y^{{\rm val(h)}};h\in\pm H(S)\}.$$
\end{thm}

Thus $\gr(D(\C[S]))$ is a semigroup ring $\C[\Sigma]$, where $\Sigma$
is a subsemigroup of $\N^2$ minimally generated by $$(d_1,0),\ldots,(d_\nu,0),
(0,d_1),\ldots,(0,d_\nu),(1,1),\{ ({\rm val}(-h),{\rm val}(h);h\in\pm H(S)\},$$

We will study study this semigroup $\Sigma$. For all this section, $S= \langle d_1, \dots, d_{\nu}\rangle $ is a numerical semigroup and $\gr(D(\C[S]))=\C[\Sigma]$.

If $z \in \Z$, denote by $\Delta_z$ the diagonal $\{(a,b) \in \N^2 ; a-b=z\}$.
 Observing the minimal set of generators of $\Sigma$, we can immediately say that all the diagonals $\Delta_s$, with $s \in \pm S$ are contained in $\Sigma$.

Since $H(S)$ is finite, $\Sigma$
is finitely generated and, setting $|H(S)|=\delta$, the number of minimal generators of $\Sigma$ is $2\nu+1+2\delta$. Thus $\Sigma$ is an affine monoid of rank 2, and
gp$(\Sigma)$, the group generated by $\Sigma$ is $\Z^2$. The normalization
of $\Sigma$, $\bar\Sigma=\{\alpha\in{\rm gp}(\Sigma);m\alpha\in\Sigma\mbox{ 
for some $m\in\N$, $m\ge1$\}}$ is $\N^2$. We denote by $\Sigma_+$ the maximal
ideal $\Sigma\setminus\{(0,0)\}$, and we set 
$T(\Sigma)=\{\tau\in{\rm gp}(\Sigma);\tau\notin\Sigma,\tau+\Sigma_+\subseteq
\Sigma_+\}$.

Observing that, if $z \in \Z$ and $({\rm val}(-z),{\rm val}(z))=(a,b)$, then 
$({\rm val}(z),
{\rm val}(-z))=(b,a)$ (cf. Lemma \ref{base} (c)), we get for $\Sigma$ 
the following symmetric property:

\begin{cor}{\rm \cite{EE},\cite{AE} }\label{simmetry}
If $a,b \in \N$, then  $(a,b)\in\Sigma$ if and only if $(b,a)\in\Sigma$.
\end{cor}

 \begin{lem}\label{generators}
 Let $S$ be a numerical semigroup.   Let $a,b \in \mathbb N$. Then 
$(a,b )\in \Sigma$ if and only if $a-b \in V_b$.
 \end{lem}

\noindent
  {\bf Proof.} We already observed that $\Sigma$ contains all the diagonals 
$\Delta_s$, for $s \in \pm S$. If $s \in S$, for each element $(a,b)\in 
\Delta_s$ (respectively  $(a,b)\in \Delta_{-s}$), we get 
val$(a-b)={\rm val}(s)=0$ and $0 \leq b$ (respectively 
val$(a-b)={\rm val}(-s)={\rm val}(s)+s=0+s=s$ and $s \leq b$). Thus in both 
cases val$(a-b) \leq b$, i.e. $a-b \in V_b$.
 
 Further, if $h \in \pm H(S)$, then 
  $({\rm val}(-h),{\rm val}(h))= ({\rm val}(h)+h,{\rm val}(h)) \in \Delta_h$ 
is a minimal generator of $\Sigma$. Since $(1,1) \in \Sigma$, it follows that  
an element $(a,b)$ of $\Delta_h$ is in $\Sigma$ if and only if 
$b \geq {\rm val}(h)={\rm val}(a-b)$, i.e. $a-b \in V_b.$ $\Box$
 
 \medskip\noindent
{\bf Definition} Let $\Gamma$ be a subsemigroup of $\N^2$, and let 
$\gamma\in\Gamma$. The Apery set of $\Gamma$ with respect to $\gamma$ is  
$\ap_\gamma(\Gamma)=\{ \alpha\in\Gamma;\alpha-\gamma\notin\Gamma\}$.

\begin{lem}\label{ap} $\ap_{(1,1)}(\Sigma)=\{(0,s);s\in S\}\cup\{(s,0);s\in S\}\cup
\{ (\val(-h),\val(h));h\in\pm H(S)\}$.
\end{lem}
 {\bf Proof.} Let $(a,b)\in\ap_{(1,1)}(\Sigma)$. We can suppose that $a\ge b$
due to the symmetry of $\Sigma$. It is clear that $(s,0)\in\ap_{(1,1)}(\Sigma)$
if and only if $s\in S$. Suppose $h\in H(S)$. We know that 
$\alpha=(h+\val(h),\val(h))\in\Delta_h$  is a minimal generator of $\Sigma$. Thus
 $\alpha\in\ap_{(1,1)}(\Sigma)$ and no other element $\sigma$ of $\Delta_h \cap \Sigma$ is in $\ap_{(1,1)}(\Sigma)$, because $\sigma -(1,1) \in \Sigma.$ $\Box$

\begin{prop}\label{T} Let $\tau\in\Z^2$. Then the following conditions are 
equivalent:

\noindent
{\rm (1)} $\tau \in T(\Sigma)$.

\noindent
{\rm (2)} $\tau+(1,1)$ is a minimal generator of $\Sigma$ of the form  
$({\rm val}(-h),{\rm val}(h))$ with $h\in\pm H(S)$.
\end{prop}

\noindent
 {\bf Proof.} (2)$\Rightarrow$(1): It is clear that $\tau\notin\Sigma$ since
$\tau+(1,1)$ is a minimal generator. It is also clear that if 
$\tau=(\tau_1,\tau_2)$, then $\tau_i\ge0$, $i=1,2$. We can suppose that
$\tau_1\ge\tau_2$ due to the symmetry of $\Sigma$. Thus $\tau=(h+i,i)$ for
some $h\in H(S)$, $i\ge0$. We have 
$\sigma=(h+\val(h),\val(h))\in\ap_{(1,1)}(\Sigma)$ according to Lemma \ref{ap}.
Let $\sigma'$ be a minimal generator. Then $\sigma+\sigma'$ is not a minimal
generator, so $\sigma+\sigma'\notin\ap_{(1,1)}(\Sigma)$, since certainly $\sigma+\sigma'\notin \{(0,s);s\in S\}\cup\{(s,0);s\in S\}$. Thus 
$\sigma+\sigma'-(1,1)\in\Sigma$, so $\tau=\sigma-(1,1)\in T(\Sigma)$.

\noindent
(1)$\Rightarrow$(2): If $\tau=(\tau_1,\tau_2)\in T(\Sigma)$, then $\tau_i\ge0$,
$i=1,2$, $\tau\notin\Sigma$, and $\tau+(1,1)\in\Sigma$. Thus 
$\tau\in\ap_{(1,1)}(\Sigma)$, so $\tau=(\val(-h),\val(h))$ for some 
$h\in\pm H(S)$ according to Lemma \ref{ap}. $\Box$

\medskip
\noindent{\bf Example} Let $S=\langle 3,5\rangle$. Then $H(S)=\{ 7,4,2,1\}$ and
$\val(7)=1,\val(4)=2,\val(2)=2,\val(1)=3$. Thus, if $\gr(D(\C[S]))=
\C[\Sigma]$, then $\Sigma$ is minimally generated by $(3,0),(5,0),(1,1),(0,3),
(0,5),(8,1),(6,2),(4,2),(4,3),(1,8),(2,6),\\
(2,4),(3,4)$. Thus 
$T(\Sigma)=\{(7,0),(5,1),(3,1),(3,2),(0,7),(1,5),(1,3),(2,3)\}$.

\medskip

Proposition \ref{T} gives a 1-1 correspondence between $T(\Sigma)$ and
the minimal generators of the form $(\val(-h),\val(h))$ with $h\in\pm H(S)$, so we get:

\begin{cor}\label{number}
The minimal set of generators of $\Sigma$ has cardinality $2\nu+1+|T(\Sigma)|$.
\end{cor}

  \medskip
 
We know by Lemma \ref{generators} that $\Sigma= \cup_{b \geq 0} (b+V_b,b)$. In a similar way, we can describe $\Sigma \cup T(\Sigma)$:

\begin{prop} \label{union}$\Sigma \cup T(\Sigma) = \cup_{b \geq 0} (b+V_{b+1},b)$.

\end{prop}

\noindent
 {\bf Proof.} If $\sigma \in \Sigma$, then, for some $ b\geq 0$, $\sigma \in (b+V_b,b)\subseteq (b+V_{b+1},b)$. If $\sigma \in T(\Sigma)$, then $\sigma +(1,1) \in  (b+V_b,b)$, for some $b \geq 1$, thus $\sigma \in ((b-1)+V_b,b-1)$, for some $b \geq 1$. Finally, if $\sigma \in \mathbb N^2 \setminus T(\Sigma)$, then $\sigma +(1,1) \notin  (b+V_b,b)$, for any $b \geq 1$, so $\sigma \notin ((b-1)+V_b,b-1)$, for any $b \geq 1.$ $\Box$
  
  \medskip
 
 \begin{prop}  Let $S$ be a numerical semigroup of maximal embedding dimension and let $\gr(D(\mathbb C[S_1])= 
\mathbb C [\Sigma_1]$. Then 

\noindent {\rm (i)} $\Sigma_1= \Sigma \cup T(\Sigma)$.

\noindent {\rm (ii)} If $a,b \in \mathbb N$, then $(a,b) \in \Sigma_1$ if and only if $(a+1,b+1) \in \Sigma$.

\end{prop}

\noindent
 {\bf Proof.} The proof follows combining  Lemma \ref{MED} and Proposition \ref{union}. $\Box$

\medskip

\noindent{\bf Example} Let $S=\langle 4,6,9,11\rangle$. Then $$T(\Sigma)=\{ (3,2),(2,0),(4,1),(5,0),(7,0),
(2,3),(0,2),(1,4),(0,5),(0,7)\}$$ 
and 
$\N^2\setminus\Sigma=T(\Sigma)\cup\{ (1,0),(3,0),(2,1),(0,1),
(0,3),(1,2)\}$. The blowup of $S$ is $S'=\langle2,5\rangle$, $\Sigma'\setminus\Sigma=T(\Sigma)$, and
$T(\Sigma')=\{(3,2),(4,1),(2,3),(1,4)\}$. The blowup of $S'$ is $S''=\langle2,3\rangle$, $\Sigma''\setminus\Sigma'=T(\Sigma')$, and $T(\Sigma'')=\{((2,1),(1,2)\}$. The blowup of $S''$ is $\N^2$,  $\N^2\setminus\Sigma''=T(\Sigma'')$.

\medskip

Now let's consider the case when the starting numerical semigroup $S$ is Arf.

\begin{lem}\label{Arfgenerators}
Let $S$ be an Arf numerical semigroup.   Let $a,b \in \mathbb N$, $a \geq b$. 
Then $(a,b)\in \Sigma$ if and only if $a-b \in I_b-s_b$.
\end{lem}
 {\bf Proof.} By Lemma \ref{generators} $(a,b)$, with $a \geq b$ is an element 
of  $\Sigma$
 when $a-b\in V_b$ and $a-b \geq 0$. By Lemma \ref{Arf} this is 
equivalent to $a-b \in I_b-I_b$. But $I_b-I_b=I_b-s_b$, because $S$ is Arf. $\Box$
 
\medskip

\begin{prop}\label{ArfT}
If $S$ is an Arf semigroup,  $S \neq \N$,  then
$$T(\Sigma)=\cup_{b=0}^{n-1}((T(S_b)+b)\times\{ b\})\cup\cup_{a=0}^{n-1}(\{ a\}
\times(T(S_a)+a)).$$
\end{prop}

\noindent
 {\bf Proof.} Let $\tau=(t+b,b)$, with $t\in T(S_b)$ and $b\in\N$. Since
$\tau\in\N^2$, to show that $\tau\in T(\Sigma)$, it is enough to show that
$\tau\notin\Sigma$ and $\tau+(1,1)\in\Sigma$. Since $t\in T(S_b)$, we have
$t\in S_{b+1}\setminus S_b$. So $\tau=(t+b,b)\in(T(S_b)+b,b)$, that is
$(t+b,b)\notin(S_b+b,b)$ and $(t+b,b)+(1,1)\in(S_{b+1},b+1)\subseteq\Sigma$.
By the symmetry of $\Sigma$ it follows that each element of the form
$(a,t+a)$, with $t\in T(S_a)$ and $a\in\N$ is in $T(\Sigma)$. For the opposite
inclusion, note that, by Lemma \ref{T}, $|T(\Sigma)|=2|H(S)|$. Counting the
elements in $\cup_{b=0}^{n-1}((T(S_b)+b)\times\{ b\})\cup\cup_{b=0}^{n-1}(\{ a\}
\times(T(S_a)+a))$, we have $2((e_0-1)+(e_1-1)+\cdots+(e_{m-1}-1)=2|H(S)|$
elements, since for an Arf semigroup $S$ with multiplicity $e$, $|T(S)|=e-1.4$ $\Box$

\begin{prop} If $S$ is an  Arf semigroup, $S \neq \mathbb N$ with multiplicity sequence $e_0,e_1,
\ldots,e_{n-1}\ne1,e_n=1,e_{n+1}=1,\ldots$,  then

\noindent
{\rm (a)} $|\N^2\setminus\Sigma|=2(e_0+2e_1+\cdots+ne_{n-1}-{n+1\choose2})
=2(ns_n-(s_1+s_2+\cdots+s_{n-1})-{n+1\choose2}) $.

\noindent
{\rm(b)} $\Sigma$ is minimally generated by $2e_0+1+2((e_0-1)+(e_1-1)+\cdots+
(e_{n-1}-1))$ elements.
\end{prop}

\noindent
 {\bf Proof.} (a) By the symmetry of $\Sigma$, 
it suffices to consider the set $\{(a,b) \in \N^2 ;a>b\}$. The number of $(a,0) \notin \Sigma $ is $s_n-(n-1)=(e_0+e_1+\cdots+e_{n-1})-(n-1)$. The number 
of $(a,1) \notin \Sigma $, $a>1$, is $e_1+e_2+\cdots+e_{n-1}-(n-2)$.
The number of $(a,2) \notin \Sigma$, $a>2$, is 
$e_2+\cdots+e_{n-1}-(n-2)$, and so on. Thus we get the left hand side. Since
$e_0+e_1+\cdots+e_{i-1}=s_i$, we have $ns_n-(s_1+s_2+\cdots+s_{n-1})=
n(e_0+e_1+\cdots+e_{n-1}-(e_0+(e_0+e_1)+\cdots+(e_0+e_1+\cdots+e_{n-2}))=
e_0+2e_1+\cdots+ne_{n-1}$.
(b) We have seen that $\Sigma$ is minimally generated by $2\nu+1+2\delta$ elements. Since $S$ is Arf, we have $\nu=e_0$ and $\delta=(e_0-1)+(e_1-1)+\cdots+(e_{n-1}-1).$ $\Box$

 \medskip
 
 \begin{prop}
If $S$ is any numerical semigroup, $S \neq \N$, the number $\mu$ of minimal generators  of $\Sigma$ satisfies
$g +6\le \mu \le 4g +3$, with equality to the left if and only if $S$ is 
2-generated 
and equality to the right if and only if 
$S=\langle g +1,g +2,\ldots,2g +1\rangle$.
\end{prop}

\noindent
 {\bf Proof.}
We know   that the number of minimal generators 
of $\Sigma$ is $2\nu +1 +2 \delta$, where $\nu$ is the number of generators for $S$,
and $\delta$ is the number of gaps. The number of gaps is at least $(g +1)/2$,
and the number of generators is at least 2. The number of gaps is at most
$g $, and the number of generators at most $ g+1$. If $\nu=2$, then $S$ is  symmetric, thus $\delta=(g +1)/2$ and we have equality to the left. On the other hand, it is  $\delta=g $ if and only if $S=\langle g +1,g+2,\ldots,2g +1\rangle$. This  is a semigroup of maximal embedding dimension of multiplicity $e=g+1$, so $\nu=e=g+1$ and we have equality to the right. $\Box$

\medskip
 
We give two examples of the ring of differential operators and its
associated graded ring for Arf semigroup rings.

\smallskip\noindent
{\bf Example 1}
If $S=\langle2,5\rangle$, then $\Sigma$ is given  by all
$(a,b)\in\N^2$ except $$\{ (1,0),(3,0),(0,1),(0,3),(2,1),(1,2)\}$$ 
A minimal set of generators for $\C[\Sigma]$ is $$\{ t^2,t^5,y^2,
y^5,ty,t^4y,t^3y^2,t^2y^3,ty^4\}$$  so $T(\Sigma)=\{(3,0),(2,1),(0,3),(1,2)\}$. 
A corresponding set of generators for the ring of differential operators $D(\C[S])$ is $$\{ t^2,t^5,\partial^2-
4t^{-1}\partial,\partial^5-10t^{-1}\partial^4+45t^{-2}\partial^3-
105t^{-3}\partial^2+105t^{-4}\partial,t\partial,$$ $$t^4\partial,
 t^3\partial^2-
t^2\partial,t^2\partial^3-3t\partial^2+3\partial,t\partial^4-6\partial^3+
15t^{-1}\partial^2-15t^{-2}\partial\}$$ This is not a minimal generating set,
e.g. $[t^4\partial,t^2]=2t^5$. The blowup of $S$ is $S_1=\langle 2,3\rangle$. All 
monomials except $t$ and $y$ belong to $\C[\Sigma_1]={\rm gr}(D(\C[S_1]))$. In general,
if $S=\langle 2,2k+1\rangle$, $|\N^2\setminus\Sigma|=k(k+1)$. 
All elements
$(a,b)$, $a,b\in \N^2$, except those where $a+b=2i-1,i=1,\ldots,k$, belong to
$\Sigma$.

\smallskip\noindent
{\bf Example 2}
If $S=\langle3,4,5\rangle$, then $\Sigma$ is given by all
$(a,b)\in\N^2$ except $$\{ (1,0),(2,0),(0,1),(0,2)\}$$ A minimal set of 
generators for $\C[\Sigma]$ is $$\{ t^3,t^4,t^5,y^3,\\
y^4,y^5,ty,t^2y,t^3y,ty^2,ty^3\}$$ so $T(\Sigma)=\{ (1,0),(2,0),(0,1),(0,2)\}$. 
A corresponding set of generators for $D(\C[S])$ is $$\{ t^3,t^4,t^5,\partial^3-
6t^{-1}\partial^2+12t^{-3}\partial,\partial^4-8t^{-1}\partial^3+
28t^{-2}\partial^2-40t^{-3}\partial,$$ $$\partial^5-10t^{-1}\partial^4+
50t^{-2}\partial^3-140t^{-3}\partial^2+180t^{-4}\partial,t\partial,
t^2\partial,t^3\partial,t\partial^2-
2\partial,t\partial^3-4\partial^2+6t^{-1}\partial\}$$ This is not a minimal 
generating set, e.g. $[t^2\partial,t^3]=3t^4$. The blowup of $S$ is $\N$ and gr$(D(\C[\N])=\C[t,y]$.

\medskip

\section{Irreducible Ideals}
It is well known that, for a numerical semigroup $S$, the cardinality of $T(S)$ is the CM type of $\C[S]$, i.e. $t=|T(S)|$ is the number of components of a decomposition of a principal ideal as irredundant intersection of irreducible ideals. We want to study whether $|T(\Sigma)|$ has a similar meaning in the ring $\C[\Sigma]$.

 Let $I$ be a proper  ideal of   $\Sigma$ i.e. a proper subset $I$ of $\Sigma$ such that $I+ \Sigma \in I$. 
 $I$ is  {\it   irreducible} if it is not the intersection of two (or, equivalently, a finite number of) ideals which properly contain $I$.
 $I$  is {\it completely   irreducible} if it is not the intersection of any set of   ideals which properly contain $I$.
 
  Consider the partial order on $\Sigma$ given by 
$$\sigma_1 \preceq \sigma_2\  \Leftrightarrow     \sigma_1+\sigma_3=\sigma_2{\rm, \ \ for\  some} \ \ \sigma_3 \in \Sigma \eqno(\star)$$
and for $x \in \Sigma$, set
$$B(x)=\{\sigma \in \Sigma \ |\ \sigma \preceq x\}$$
  
\begin{lem} \label{cirred} If $I$ is a proper  ideal of   $\Sigma$, then the following conditions are equivalent:
 
 {\rm (1)} $I$  is completely irreducible. 

{\rm (2)} $I$ is maximal as  ideal with respect to the property of not containing an element $x$, for some $x \in \Sigma$.

{\rm (3)} $I=\Sigma \setminus B(x)$, for some $x \in \Sigma$.
 \end{lem}
 
 \noindent
 {\bf Proof.}   $(1)\Rightarrow (2)$: Let $H$ be the intersection of all the  ideals properly containing $I$. Then there is $x \in H \setminus I$, so $I$ is maximal with respect to the property of not containing $x$.  
 
 \noindent
 $(2)\Rightarrow (1)$. Each  ideal $J$ properly containing $I$ contains $x$, so $I$ is not the intersection of all such ideals $J$ and it is completely irreducible. $(2) \Leftrightarrow (3)$ is trivial.  $\Box$

\begin{lem}\label{irred} For each $a,b \in \mathbb N$, $a,b>0$, the following are irreducible, non completely irreducible ideals of $\Sigma$:

$N_{(a,0)}:= \Sigma \cap\{{(x,y)\in \N^2;x\ge a}\}$

$N_{(0,b)}:= \Sigma \cap\{{(x,y)\in \N^2;y\ge b}\}$
\end{lem}

\noindent
{\bf Proof.}  Any ideal $J$ of $\Sigma$ properly containing $N_{(a,0)}$ contains $(a-1,s)$, for some $s \in S$, so it contains $(a-1,s+S)$. It follows that, if $J_1$, $J_2$ are ideals properly containing $N_{(a,0)}$, then $J_1 \cap J_2$ contains $(a-1, \max(s,s')+g+1+S)$ (if $(a-1,s) \in J_1$ and $(a-1,s') \in J_2$), so $J_1 \cap J_2 \neq N_{(a,0)}$ and $N_{(a,0)}$ is irreducible. On the other hand $N_{(a,0)}$ is the intersection of all the ideals properly containing it, so it is not completely irreducible. $\Box$

\medskip

In all the results of this section $\max$ has to be intended with respect to the partial order $(\star)$ on $\Sigma$.

\begin{prop}\label{intersection} Let $I$ be an ideal of $\Sigma$ generated by $(a_1,b_1), \dots, (a_h,b_h)$ and let $a= \min\{a_i\}$, $b= \min\{b_i\}$. Then
$$ I= \bigcap_{x \in \max (\Sigma \setminus I)} (\Sigma \setminus B(x)) \cap N_{(a,0)} \cap N_{(0,b)}$$
is the unique irredundant decomposition of the ideal $I$ as intersection of irreducible ideals.
\end{prop}
{\bf Proof.} $\subseteq$: let $\alpha \in I$. Then $\alpha \notin B(x)$ for each $x \in \Sigma \setminus I$ (otherwise $\alpha + \beta =x$, for some $\beta \in \Sigma$ and so $x \in I$, a contradiction). Thus $\alpha \in \Sigma \setminus B(x)$, for each $x \in \max (\Sigma \setminus I)$. Moreover $\alpha \in N_{(a,0)} \cap N_{(0,b)}$.
$\supseteq$: observe first that $$\bigcap_{x \in \max (\Sigma \setminus I)}(\Sigma \setminus B(x))=\bigcap_{x \in  (\Sigma \setminus I)}(\Sigma \setminus B(x))$$
in fact $(\Sigma \setminus B(x_1)) \subseteq (\Sigma \setminus B(x_2))$ if and only if $B(x_1) \supseteq B(x_2)$ if and only if $x_2 \preceq x_1$.
Suppose that $\alpha \in  N_{(a,0)} \cap N_{(0,b)}$, i.e. that $\alpha=(c,d) \in \Sigma$, with $c \ge a$ and $d \ge b$. We have to show that, if $\alpha \in \bigcap_{x \in  (\Sigma \setminus I)}(\Sigma \setminus B(x))$, then $\alpha \in I$. In fact, if $\alpha \notin I$, then (since trivially $\alpha \in B(\alpha)$) $\alpha \notin \Sigma \setminus B(x)$, for some $x \in \Sigma \setminus I$ (take $x=\alpha$).

To show that the decomposition is irredundant, it's easy to see that $N_{(a,0)}$ (respectively $N_{(0,b)}$) does not contain the intersection of the other components. Moreover, if $x \in \max (\Sigma \setminus I)$, the only component of the intersection which does not contain $x$ is $\Sigma \setminus B(x)$. Thus this component is not superfluos. $\Box$

\medskip

The following result agrees with \cite[Theorem 11.3]{MS}:
\begin{cor} The unique irreducible ideals of $\Sigma$ are $N_{(a,0)}$, for some $a>0$, $N_{(0,b)}$, for some $b>0$ and those of the form $\Sigma \setminus B(x)$, for some $x\in \Sigma$.

\end{cor}

\begin{cor} If $(0,0) \neq \sigma=(a,b) \in \Sigma$, then 
$$ \sigma + \Sigma = \bigcap_{x \in \max {\rm Ap}_{\sigma}(\Sigma)} (\Sigma \setminus B(x)) \cap N_{(a,0)} \cap N_{(0,b)}$$
is the unique irredundant decomposition of the principal ideal $\sigma + \Sigma$ as intersection of irreducible ideals.

\end{cor}
{\bf Proof.} It follows from the Proposition \ref{intersection}, observing that ${\rm Ap}_{\sigma}(\Sigma)= \Sigma \setminus (\sigma+ \Sigma).$ $\Box$

\medskip

What we got for the ideals of the semigroup $\Sigma$  can be read in terms of monomial ideals of $\C[\Sigma]$. In fact each ideal $I$ of $\Sigma$ corresponds to the monomial ideal of $\C[\Sigma]$ generated by $\{t^ay^b;\ (a,b)\in I\}$. Moreover if a monomial ideal of $\C[\Sigma]$ is not the intersection of two strictly larger monomial ideals, then it is not the intersection of two strictly larger ideals, even if non monomial ideals are allowed (\cite[Proposition 11, p.41]{MS}). Thus the results above characterize the irreducible monomial ideals of $\C[\Sigma]$ as well.

\begin{cor} Each principal monomial ideal of $\C[\Sigma]$ is an irredundant intersection of $|T(\Sigma)|+2=2 \delta+2$ irreducible ideals (where $\delta =|H(S)|)$.

\end{cor}
{\bf Proof.} It follows from  the previous Corollary, recalling that, for each $(0,0) \neq \sigma \in \Sigma$,  $2 \delta=|T(\Sigma)|=| \max  {\rm Ap}_{\sigma} (\Sigma)|$, because there is a one to one correspondence between the sets $T(\Sigma)$ and $ \max  {\rm Ap}_{\sigma} (\Sigma)$, more precisely it is proved in \cite[Proposition 4.1]{E} that $\tau \in T(\Sigma)$ if and only if $\tau + \sigma \in \max {\rm Ap}_{\sigma} (\Sigma).$ $\Box$

\section{Derivations}
Let $I$ be an ideal in $\C[S]$. Then we denote by $\der(I,I)$ the set of 
derivations which map $I$ into $I$.
\begin{lem}\label{der}
If $I$ is generated by monomials in $\C[S]$, then $\der(I,I)\simeq I:I$ as
$\C[S]$-module. Thus $\der(I,I)$ is isomorphic to a semigroup ring $\C[T]$, where $T$ is a semigroup, 
$S\subseteq T\subseteq\N$.
\end{lem}

\noindent
{\bf Proof.}  If $I$ is generated by monomials 
also $I:I$, which is a fractional ideal of $\R$, is generated by 
monomials.
Let
$\{ t^{n_i}\}$ be the generators of $I:I$, then $\der(I,I)$ is
generated by $\{ t^{n_i+1}\partial\}$, and $t^k\mapsto t^{k+1}\partial$
induces an isomorphism as $\C[S]$-modules.  Moreover $I:I$ is a semigroup ring $\C[T]$, for some semigroup $T$, 
$S\subseteq T\subseteq\N$ and so  $\der(I,I)$ is isomorphic to   $\C[T]$ as $\C[S]$-module.
$\Box$

\smallskip
Observe that, if $I:I=\C[T]$, then also $xI:xI=\C[T]$, for each nonzero $x\in \R$.
In particular we can say that, for each monomial principal ideal $I$ of $\C[S]$, 
  $\der(I,I)\simeq \C[S]$ as 
$\C[S]$-modules.

 \smallskip
If $I$ is not generated by monomials, the statement in the proposition is no 
longer true. If $I=(t^4+t^5,t^4+t^6)$ in $k[t^4,t^5,t^6]$, then $\der(I,I)$ is 
generated by $t^5\partial,t^6\partial,t^7\partial$.

\bigskip

For a monomial ideal $I$ of $\C[S]$, we denote by $\min(I)$ the minimal degree of the monomials in $I$.  

\begin{prop}\label{min}
Let $T$ be a semigroup, $S \subseteq T  \subseteq \N$ and let $C$ be the conductor ideal $C= \C[S]:\C[T]$. Then:

\noindent {\rm (i)} If $J$ is a monomial ideal of $\C[S]$ such that $J:J= \C[T]$, then $J \subseteq C$ and $\min(J) \geq \min(C)$.

\noindent {\rm (ii)} The overring $\C[T]$ of $\C[S]$ is of the form $I:I$, for some monomial  ideal $I$ of $\C[S]$.

\noindent {\rm (iii)} If $J$ is a monomial ideal of $\R$, with $\min(J)=j$, such that $J:J=\A$, then $J \supseteq t^j\A$ and $t^j \in C$.

 \end{prop}

\noindent
{\bf Proof.} (i)  If $J$ is a monomial ideal of $\C[S]$ such that $J:J= \C[T]$, then, since $J$ is an ideal of $\A$ too, it is contained in $C$ which is the biggest ideal that  $\R$ and $\A$ share. Thus $\min(J) \geq \min(C)$.

(ii)  We have that $\A \subseteq \overline{\C[S]} = \C[t]=\C[\N]$, so $\A$ is a fractional ideal of $\R$ and $r\A=I \subseteq \R$, for some nonzero  $r \in \R$. Now $I:I= r\A:r\A =\A$ and $\A$ is of the requested form.

(iii) Suppose now that $J$ is  a monomial ideal of $\R$, with $\min(J)=j$, such that  $J:J=\A$.  We claim that $J \supseteq t^j\A$. Indeed the principal ideal  $t^j \R$ is contained in $J$ and so $\A=J:J \subseteq J: t^j\R= t^{-j}(J: \R)=t^{-j}J$, hence $ t^j\A \subseteq J$. Then also $ t^j\A \subseteq \R$ and $t^j \in C$. $\Box$

\medskip
It is known that  the ring $\C[[T]]$ is Gorenstein if and only if the numerical semigroup $T$
 is symmetric. The extension of that result to the non local case is not difficult: 

 \begin{lem}\label{Kunz} Let $T$ be a numerical semigroup. Then the ring $\C[T]$ is Gorenstein if and only if $T$ is symmetric.
\end{lem}

{\bf Proof.} $\C[T]=\C[t^{n_1}, \dots, t^{n_h}]$ is
 Gorenstein if and only if each localization at a prime ideal is Gorenstein. For the
 localization at $ P=( t^{n_1}, \dots, t^{n_h})$, we have that $\C[T]_P$ is Gorenstein
 if and only if $T$ is symmetric, argueing similarly to the local case $\C[[T]]$. For the other nonzero
 prime ideals $Q$, we have that $Q$ does not contain the conductor $\C[T]:\C[t]$,
so $\C[T]_Q \cong \C[t]_{Q'}$ is a DVR, thus a Gorenstein ring, where $Q'$ is the unique prime ideal of
$\C[t]$ lying over $Q$. $\Box$

\begin{cor}\label{Gorenstein}
Let  $\A$ be a Gorenstein overring of $\R$. If $J$ is a monomial ideal of $\R$, then $J:J= \A$ if and only if $J=t^j\A$, where $j=\min(J)$.

\end{cor}

\noindent
{\bf Proof.}  We know by Proposition \ref {min} iii) that $J \supseteq t^j\A$. In order to prove the opposite inclusion, we show that, if $H$ is a monomial ideal of $\R$ with $\min(H)=j$ properly larger than $t^j\A$, then  $H:H \neq \A$. We can argue equivalently on the fractional ideal $\A$ of $\R$. Let $g$ be the Frobenius number of the semigroup $T$, which is symmetric because $\C[T]$ is Gorenstein. So, let $H$ be a fractional monomial ideal of $\R$ strictly larger than $\A$ and let $t^h \in H \setminus \A$. Then $t^{g-h} \in \A$ and so $t^{g-h} t^h=t^g \notin \A$. Thus $t^{g-h} \notin H:H$ and $H:H \neq \A$. $\Box$
 \smallskip

  If $I$, $J$ are ideals of $\R$, we say that $I$ and $J$ are {\it equivalent} if $xI=yJ$, for some nonzero elements $x,y \in \R$.
Recall  also that an ideal $I$ is called {\it stable} if it is principal in the overring $I:I$.  

\smallskip\noindent
{\bf Example} If $S=\langle 3,4,5\rangle$, the only  semigroups $T$ for which $S\subseteq T\subseteq\N$ are
$S$, $T_1=\langle 2,3\rangle$, and $\N$, so that the only proper semigroup overrings of $\R$ are $\C[T_1]$ and $\C[\N]=\C[t]$, which are both Gorenstein rings. 
The conductor ideals coincide, $\C[S]:\C[T_1]=\C[S]:\C[\N]=t^3\C[t]$. By Corollary \ref{Gorenstein}, if $J$ is a monomial ideal of $\C[S]$, then $J:J=\C[T_1]$ if and only if $J=t^j\C[T_1]$ with
$j\ge3$, i.e. $J=(t^3,t^5)t^k$, with $k \geq 0$ and $J:J =\C[t]$ if and only if $J=t^j\C[t]$ with $j\ge3$, i.e. $J=(t^3,t^4,t^5)t^k$, with $k \geq 0$. On the other hand, the trivial overring $\R$ is not Gorenstein and $\R:\R=\R$. We have that $J:J=\R$ if and only if $J$ is a principal ideal of $\R$ or $J=(t^3,t^4)t^k$, with $k \geq 0$,  so here two equivalence classes of monomial ideals correspond to the same overring.

 \begin{prop} The following conditions are equivalent:

\noindent {\rm (1)}There exists a one-to-one correspondence between the semigroup overrings of $\R$ and the equivalence classes of monomial ideals of $\R$.

\noindent {\rm (2)} $\R= \C[t^2, t^{2k+1}]$, for some $k \in \mathbb N$.

\noindent {\rm (3)} Each semigroup overring of $\R$ is Gorenstein.
\end{prop} 

\noindent
{\bf Proof.}  (1) $\Leftrightarrow$ (2). By Proposition \ref{min} (iii), there exists a one-to-one correspondence between the semigroup overrings of $\R$ and the classes of stable monomial ideals of $\R$ (cf. also \cite [Proposition II.4.3]{bdf}). So we get the requested one-to-one correspondence if and only if each monomial ideal of $\R$ is stable, i.e. if and only if each semigroup ideal of $S$ is stable. By \cite[Theorems I.5.13, (i) $\Leftrightarrow$ (iii) and I.4.2 (i) $\Leftrightarrow$ (v)]{bdf}, this is equivalent to $S= \langle 2, 2k+1\rangle$, for some $k \in \mathbb N$. (2) $\Leftrightarrow$ (3). By Lemma \ref{Kunz}, condition (3) means that each semigroup $T$, $S \subseteq T \subseteq N$ is symmetric and that holds if and only if $S= \langle 2, 2k+1\rangle$, for some $k \in \mathbb N$  (cf. \cite [Theorem I.4.2 (v) $\Leftrightarrow$ (ix)]{bdf}). $\Box$

\medskip

\smallskip\noindent
{\bf Example.}  
 If $S=\langle2,5\rangle$ there are three equivalence classes of ideals
generated by monomials with representatives $\R$, $(t^2,t^5)$ and $(t^4,t^5)$.
These correspond to $\R$, $\C[t^2,t^3]$ and $\C[t]$, respectively. More generally, if $S=\langle 2,  2k+1\rangle$, there are $k+1$ equivalence classes of ideals
generated by monomials with representatives  $\R$, $(t^2,t^{2k+1}),\ (t^4,t^{2k+1}), \dots, \  (t^{2k},t^{2k+1})$. These correspond to $\R,\ \C[t^2, t^{2k-1}], \  \C[t^2, t^{2k-3}], \dots, \ \C[t^2, t ]=\C[t]$, respectively.

\end{document}